%% file: VennDiagramsAreHamiltonianArXiv.tex
\documentclass[12pt]{llncs}
\usepackage{llncsdoc}

\usepackage{enumitem}
\usepackage{makeidx}  
\usepackage{amsfonts}
\usepackage{xcolor}
\usepackage{graphicx}
\usepackage{lipsum}
\usepackage{verbatim}

\usepackage{tikz}
\usetikzlibrary{backgrounds}
\usetikzlibrary{intersections}

\definecolor{myRed}{RGB}{230,20,37}
\definecolor{myGreen}{RGB}{19,139,74}
\definecolor{myOrange}{RGB}{239,103,27}
\definecolor{myBlue}{RGB}{21,69,153}

\begin{document}
\frontmatter          
\pagestyle{headings}  
\addtocmark{All Simple Venn Diagrams are Hamiltonian} 

\mainmatter              
\title{All Simple Venn Diagrams are Hamiltonian } 
\titlerunning{All Simple Venn Diagrams are Hamiltonian}  
%
\author{Gara Pruesse\inst{1} \and Frank Ruskey\inst{2}}
\authorrunning{Pruesse and Ruskey} 
%
\tocauthor{Gara Pruesse, Frank Ruskey}
\institute{Vancouver Island University,\\
Department of Computing Science,\\Nanaimo, BC, Canada,\\
\email{Gara.Pruesse@viu.ca}\\
\and
University of Victoria,\\
Department of Computer Science,\\
Victoria, BC, Canada}


\maketitle              

\begin{abstract}
An $n$-Venn diagram is a certain collection of $n$ simple closed
curves in the plane.  They can be regarded as graphs where the
points of intersection are vertices and the curve segments between
points of intersection are edges.
Every $n$-Venn diagram has the property that
a curve touches any given face at most once between the
points of intersection incident to that face.
We prove that any connected collection of $n$ simple closed curves satisfying
that property are 4-connected, if $n \ge 3$, so long as the curves
intersect transversally and at most two curves intersect at any point.  Hence by a theorem of Tutte,
such collections, including simple Venn diagrams, are Hamiltonian.
\end{abstract}


\pagestyle{myheadings}

\thispagestyle{plain}
\markboth{PRUESSE AND RUSKEY}{Venn Diagrams are Hamiltonian}

\section{Introduction}
Perhaps the most intriguing open question about Venn diagrams is
that of Winkler \cite{Winkler1984}: Is every simple $n$-Venn diagram
\emph{extendible} to a simple $(n+1)$-Venn diagram by the addition of a
suitable curve?  Winkler conjectured (``nervously") that the answer is
yes.

An $n$-Venn diagram is a collection of such curves with the further properties that
the collection partitions the plane into $2^n$ connected non-empty regions,
each region being formed by one of the $2^n$ possible intersections of the
interior or exterior of each curve.  A \emph{simple Venn diagram} is one in
which there are a finite number of intersection points and at most two curves
pass through any given point.  In this paper we will assume that all Venn diagrams
are simple.  The figure below shows a $4$-Venn diagram on
the right.

\includegraphics[width=4in]{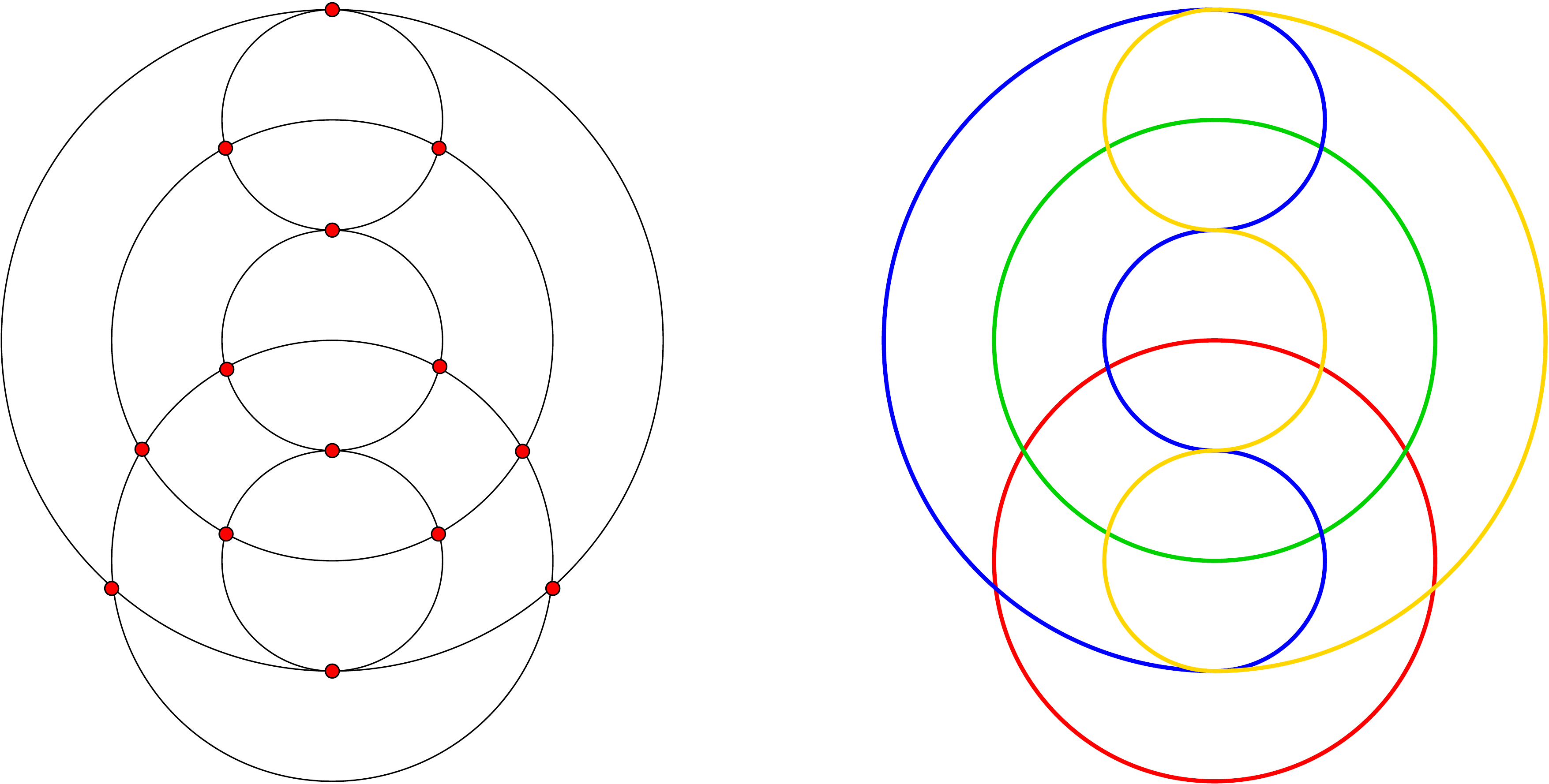}

A Venn diagram can be regarded as a planar graph in which the vertices are
the points of intersection and edges are the segments of the curves
between the vertices; we call such a graph a \emph{Venn graph}.
The Venn graph corresponding to the Venn diagram to the right
is shown on the left in the figure above.

In the Venn graph context, Winkler's question is equivalent to
asking whether the planar dual of every simple Venn diagram is Hamiltonian.
A Hamilton cycle in the dual becomes the new curve.

It is also natural to ask whether every Venn graph is
Hamiltonian.  In this paper we introduce a much wider
class of configurations of simple closed curves in the plane and
prove that the graphs that arise from these configurations have
connectivity of four.  A theorem of Tutte \cite{Tutte1956} then implies
that they, and hence Venn graphs, are Hamiltonian.

This paper is concerned with collections of simple closed curves 
in the plane.  As is well known, any simple closed curve partitions the plane into
the curve itself, the \emph{interior} of the curve, and the \emph{exterior}
of the curve.  We will assume that the curves lie in \emph{general position},
meaning that they intersect at some finite
number of points, at most two curves intersect at any point,
 and that each intersection is transverse; that is, the
curves cross each other when they intersect, no ``kissing" allowed.
A \emph{2-face} in an embedding of closed curves is a face that is incident to exactly two curves.

With any collection of closed curves in general position, there is a natural way to associate a 4-regular plane graph
whose vertices are the points of intersection and whose edges are the
curve segments between the points of intersection.
For a vertex $v$ on a curve $C$, we may refer to $v$'s ``$C$ neighbours'', meaning those two neighbours of $v$ that are also on $C$.
Observe that the curves can be recovered from the graph since the intersections are transverse.
In general it is difficult to determine if the graph is Hamiltonian since
Iwamoto and Toussaint have shown that the problem is NP-complete \cite{IT94}.

In the case where the Venn diagrams are not required to be simple,
Chilakamarri, Hamburger, and Pippert \cite{CHP96a} showed them all to be extendible to a non-simple Venn diagram.
These same authors \cite{CHP96b} also showed that the Venn graph is always 3-connected and
similarly the Venn dual is always 3-connected (if there are more than 2 curves).
We strengthen this result by showing that the connectivity of a Venn graph is precisely 4.
Bultena \cite{Bultena} verified that Winkler's conjecture is true for all Venn diagrams with 5 or fewer curves.

\section{Results}

Two notable properties of simple $n$-Venn graphs for $n > 2$ are the following:
(a) there are no 2-faces; (b) a curve is incident on a face at most once.
In fact, for a general connected collection of curves, (b) implies (a) if $n > 2$ (draw a picture to
convince yourself).  We will call property (b) the \emph{UFI} property, an acronym for
unique face incidence property.

\begin{definition}
A \emph{V-graph} is the graph of a collection of at least three topologically connected curves in the plane, lying in general position,
and which satisfies the UFI property.
\end{definition}

Our aim is to prove that every V-graph is 4-connected and then invoke a theorem of Tutte \cite{Tutte1956} to infer that
V-graphs are Hamiltonian.

\begin{theorem}[Tutte]
Any 4-connected planar graph is Hamiltonian.
\end{theorem}

We make use of Menger's characterization of $k$-connectivity: A graph is $k$-connected if and only if between every pair of
distinct vertices there are at least $k$ disjoint paths.  (Disjoint here means that they have no intermediate
vertices in common, although they share endpoints.)

\begin{lemma}
\label{lemma:k-conn}
A connected graph is $k$-connected if and only if for all pairs of vertices $u,v$ at distance two from one another, there exist $k$ vertex-disjoint $u,v$-paths.
\end{lemma}

\medskip
\noindent
{\bf Proof:}
The ``only if" direction is obvious, so we consider the claim in the ``if" direction.
Suppose $G$ is graph, and all pairs of vertices at distance two in $G$ are connected by
$k$ disjoint paths.  Let $X$ be a cut of minimum size in $G$, and let $x$ be a vertex in $X$.
Thus $G\setminus (X \setminus \{x\})$ is connected, but the further removal of $x$
does disconnect the graph.  Therefore $x$ has a neighbour in each connected component of $G\setminus X$; let $u$ and $v$
be two such neighbours in different components.  Hence $u$ and $v$ are at distance two from one another in $G$, and so have at least $k$ disjoint paths between
them.  Each of these disjoint paths contains at least one member of $X$, so $|X| \ge k$. \hfill $ \Box$

\begin{theorem}
The connectivity of a V-graph is four.
\label{thm:Vconn4}
\end{theorem}

\medskip
\noindent
{\bf Proof:}
Since a V-graph is 4-regular the connectivity can not be higher than four.
We will show that the connectivity is at least 4 by using Lemma \ref{lemma:k-conn}.
Let $G$ be a V-graph, and let $u$ and $v$ be any two vertices at distance two in $G$, with $z$ being a vertex
that is a neighbour of both $u$ and $v$.

\begin{figure}
\scalebox{0.25}{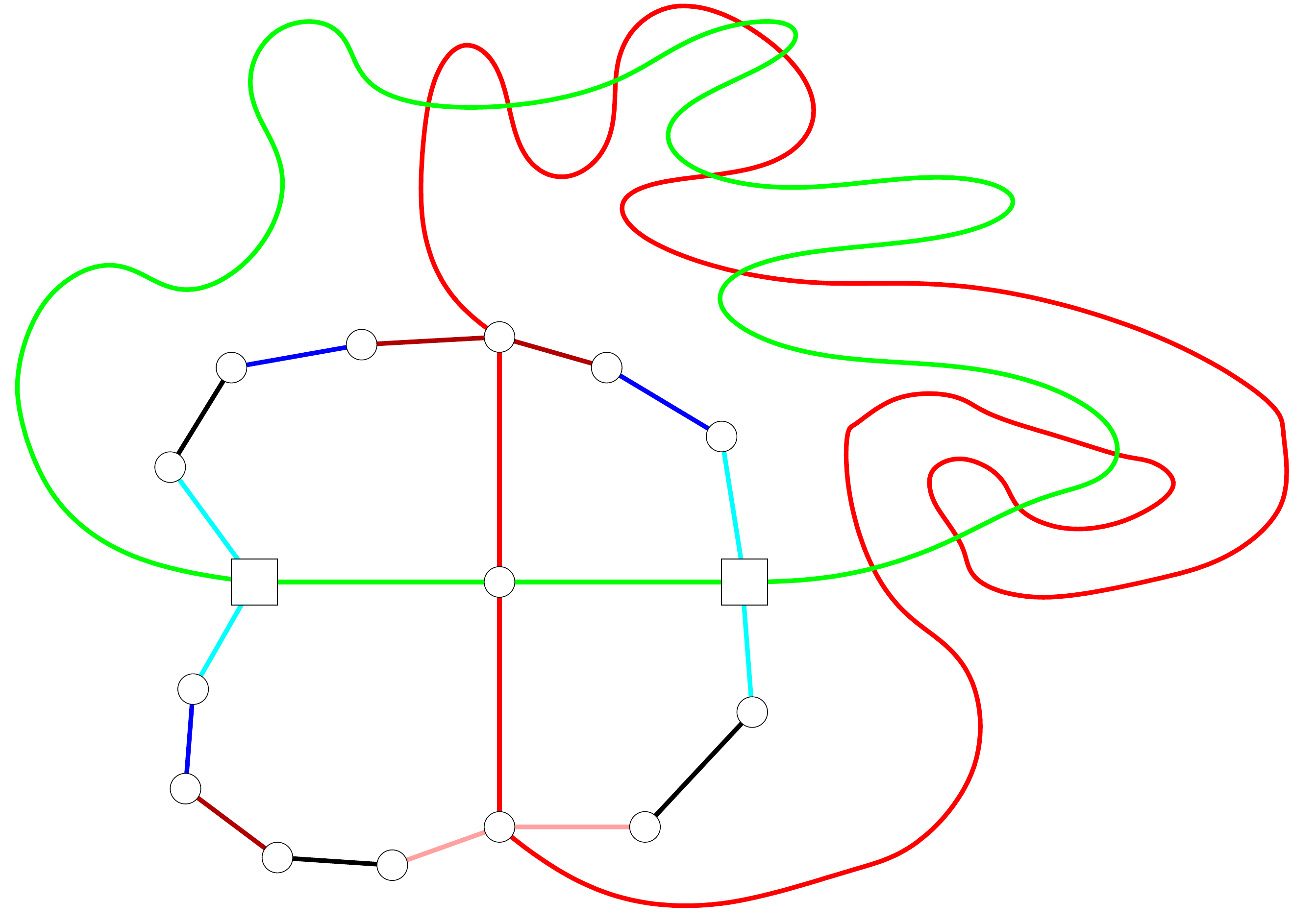}%
\scalebox{0.25}{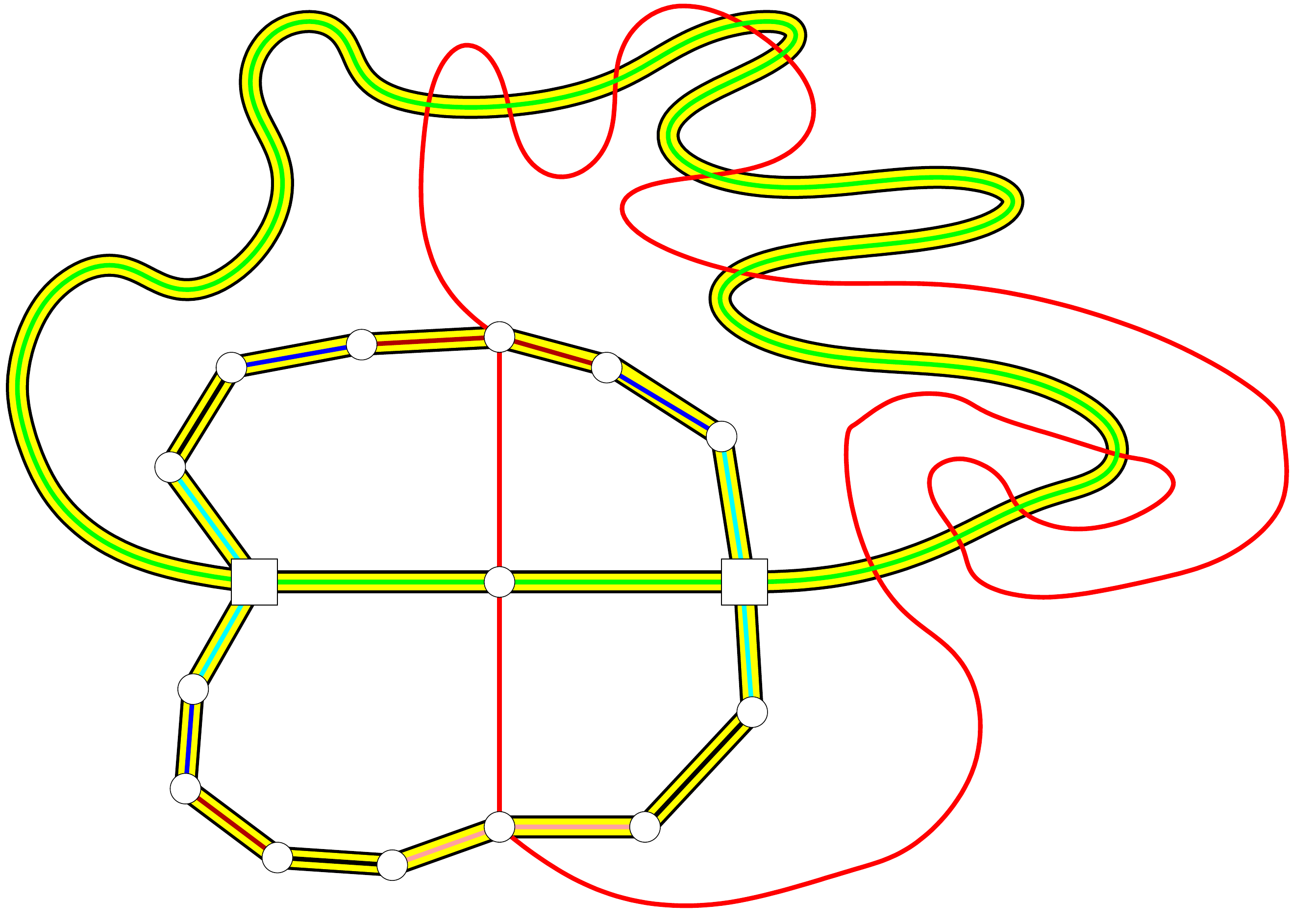}

\hfill (a) \hfill (b)\hfill
\caption{Case 1: (a) $u$ and $v$ are on the same curve. (b) Four vertex disjoint paths.}
\label{fig:case1}
\end{figure}

\begin{figure}
\scalebox{0.25}{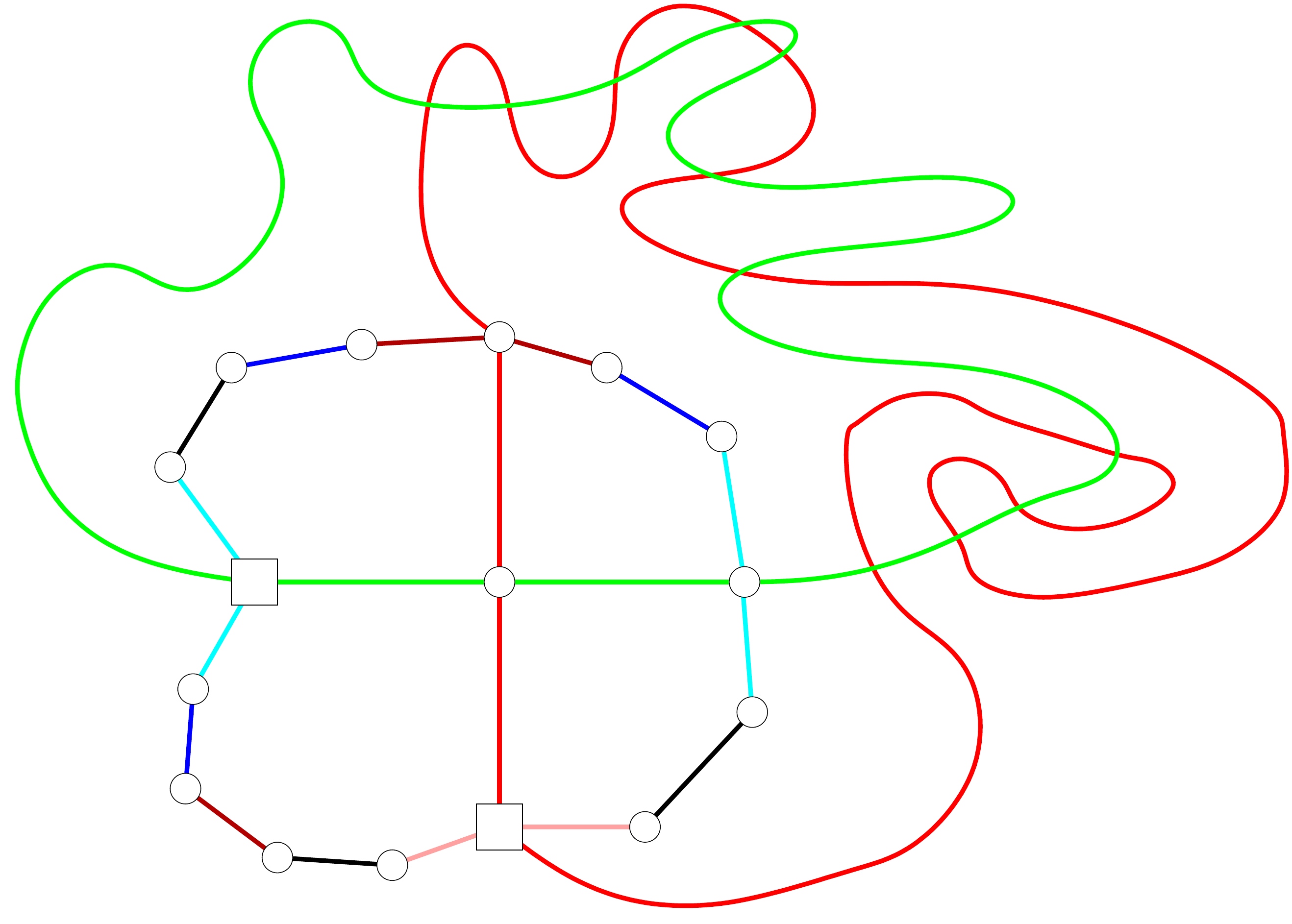}%
\scalebox{0.25}{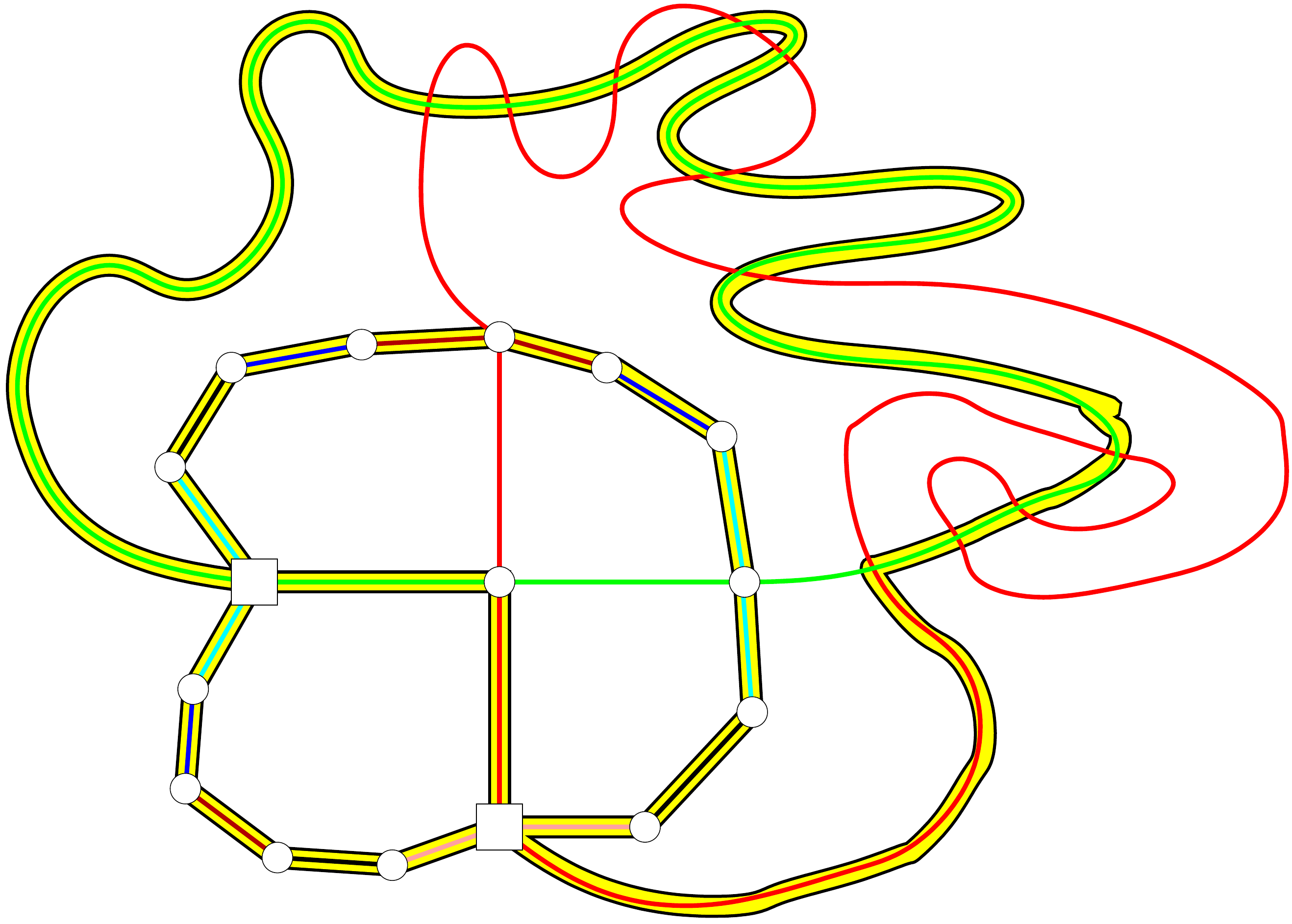}

\hfill (a) \hfill (b)\hfill
\caption{Case 2: (a) $u$ and $v$ are on different curves. (b) Four vertex disjoint paths.}
\label{fig:case2}
\end{figure}

Vertex $z$ is the intersection of exactly two curves, call them red and green.  Without loss of generality, we consider the two cases, either $u$ and $v$
are both on the green curve (Case 1),  or alternatively, $u$ is on the green curve and $v$ is on the red curve (Case 2).  We show in either case that there
are four disjoint $u,v$-paths in $G$.

\medskip
\noindent
\textbf{Case 1:} $u$ and $v$ are on the green curve.  Then $z$ has two remaining neighbours, call them $a$ and $b$, on the red curve.
Since there are no 2-faces $a$ and $b$ are distinct from $u$ and $v$ (and from each other).  See Figure \ref{fig:case1}(a).
We construct two vertex-disjoint $u,v$-paths in $G\setminus \{z\}$ by traversing the edges along the perimeter of the four
faces adjacent to $z$: the $u \rightsquigarrow a \rightsquigarrow v$ path on one side of the green curve, and the $u \rightsquigarrow b \rightsquigarrow v$
path on the other side of the green curve.   Since these paths are on different sides of the green curve they are therefore vertex disjoint from
each other, except in the endpoints.  Furthermore, by the UFI property, none of the edges on those paths are green, nor are any of the
intermediate vertices incident to a green edge.
Therefore, there are two $u,v$-paths that are disjoint from those outlined above: the two paths constructed from the cycle of green edges.
One path is $u \rightarrow z \rightarrow v$ and the other consists of the remaining edges along the green curve.
See Figure \ref{fig:case1}(b).

\medskip
\noindent
\textbf{Case 2:} $u$ is on the green curve and $v$ is on the red curve.
Let $a$ be $z$'s other red neighbour, and let $b$ be $z$'s other green neighbour.
Since there are no 2-faces, all four of $a$, $b$, $u$, $v$ are distinct.
See Figure \ref{fig:case2}(a).

We construct two vertex-disjoint $u,v$-paths in $G\setminus \{z\}$ by traversing the edges along the perimeter of the four
faces adjacent to $z$.   Call this perimeter the \emph{$z$-face} in $G\setminus \{z\}$.  Note, as before, that none of those edges is red or green.

\medskip
\noindent
\emph{Path A:} Use the non-red/green edges in the face of $G$ containing $u$ and $v$.  Denote this path $u \rightsquigarrow v$.

\noindent
\emph{Path B:} Use the remaining edges in the $z$-face: $u \rightsquigarrow a \rightsquigarrow b \rightsquigarrow v$.

\medskip
These two paths are vertex-disjoint except in their endpoints, as each non-$u,v$ vertex on Path A is separated from any non-$u,v$ vertex on Path B by either the red curve or
the green curve or both.

\medskip
\noindent
\emph{Path C:} Follow the green curve from $u$ in the direction away from $z$.  Since the green and red curves intersect at $z$ they
must intersect again; follow the green curve until the last red-green intersection.  Then follow the red curve
towards $z$ until it hits $v$.

\noindent
\emph{Path D:} Use the two edge path $u \rightarrow z \rightarrow v$.

\medskip
Paths C and D are disjoint, as they contain disjoint parts of the green and red curves.  They are disjoint from paths A and
B as they visit only vertices on the green and red curves, while avoiding vertices $a$, $b$, and $z$.

In either of Cases 1 or 2 above, there are four vertex-disjoint paths connecting $u$ and $v$, an arbitrary pair of vertices
at distance two from each other in $G$.  By Lemma 1, this suffices to prove that $G$ is 4-connected. \hfill $\Box$

\begin{theorem}
All V-graphs are Hamiltonian.
\end{theorem}

\medskip
\noindent
{\bf Proof:}
From Theorem \ref{thm:Vconn4}, above, we have that V-graph is 4-connected; and by definition, V-graphs are
planar.  By Tutte's Theorem, every 4-connected planar
graph is Hamiltonian.  Therefore every V-graph is Hamiltonian.
\hfill $\Box$

\medskip
\noindent
Since V-graphs generalize Venn graphs we have the following:

\begin{corollary}
The connectivity of a Venn graph is four.
Furthermore, every Venn graph is Hamiltonian.
\end{corollary}

\smallskip
\noindent

\section{Final Remarks}

The question of Hamiltonicity remains open if the assumption that the Venn diagrams are simple is lifted.
We know of no non-simple Venn diagrams that are not Hamiltonian.  Thus in the spirit of
Winkler's conjecture, we (nervously) make the following conjecture.

\begin{conjecture}
All non-simple Venn diagrams are Hamiltonian.
\end{conjecture}

Algorithmically, it is worth noting that there is a linear time algorithm for finding a Hamilton cycle
in a 4-connected graph; see Chiba and Nishizeki \cite{ChibaNishizeki}.
Since the problem of determining whether a collection of Jordan curves in the plane in general position
is Hamiltonian is NP-complete \cite{IT94}, it is interesting that the UFI property is enough to guarantee Hamiltonicity.

An easy example of a collection of curves that does not satisfy the UFI property
and which are only biconnected is obtained by having two curves ``weave" back-and-forth
across each other.  Such configurations have 2-faces.
Figure 7(c) on page 186 in \cite{IT94} shows an example of a non-Hamiltonian 3-connected configuration
of curves in general position that is without 2-faces.

\section{Acknowledgements}
The authors wish to thank Moshe Rosenfeld for helpful discussions and the reference to the relevant paper of Tutte \cite{Tutte1956}.
The research of the second author is supported by a Discovery Grant from NSERC.

%
%

\end{document}

%% file: HamVenn1a.pdf_t
\begin{picture}(0,0)%
\includegraphics{HamVenn1a.pdf}%
\end{picture}%
\setlength{\unitlength}{3947sp}%
\begingroup\makeatletter\ifx\SetFigFont\undefined%
\gdef\SetFigFont#1#2#3#4#5{%
  \reset@font\fontsize{#1}{#2pt}%
  \fontfamily{#3}\fontseries{#4}\fontshape{#5}%
  \selectfont}%
\fi\endgroup%
\begin{picture}(12645,9036)(4709,-9110)
\put(9601,-8911){\makebox(0,0)[b]{\smash{{\SetFigFont{34}{40.8}{\familydefault}{\mddefault}{\updefault}{\color[rgb]{0,0,0}$b$}%
}}}}
\put(11551,-6361){\makebox(0,0)[b]{\smash{{\SetFigFont{34}{40.8}{\familydefault}{\mddefault}{\updefault}{\color[rgb]{0,0,0}$v$}%
}}}}
\put(9826,-3136){\makebox(0,0)[b]{\smash{{\SetFigFont{34}{40.8}{\familydefault}{\mddefault}{\updefault}{\color[rgb]{0,0,0}$a$}%
}}}}
\put(7426,-6436){\makebox(0,0)[b]{\smash{{\SetFigFont{34}{40.8}{\familydefault}{\mddefault}{\updefault}{\color[rgb]{0,0,0}$u$}%
}}}}
\put(9301,-5536){\makebox(0,0)[b]{\smash{{\SetFigFont{34}{40.8}{\familydefault}{\mddefault}{\updefault}{\color[rgb]{0,0,0}$z$}%
}}}}
\end{picture}%

%% file: HamVenn1b.pdf_t
\begin{picture}(0,0)%
\includegraphics{HamVenn1b.pdf}%
\end{picture}%
\setlength{\unitlength}{3947sp}%
\begingroup\makeatletter\ifx\SetFigFont\undefined%
\gdef\SetFigFont#1#2#3#4#5{%
  \reset@font\fontsize{#1}{#2pt}%
  \fontfamily{#3}\fontseries{#4}\fontshape{#5}%
  \selectfont}%
\fi\endgroup%
\begin{picture}(12645,9036)(4709,-9110)
\put(9601,-8911){\makebox(0,0)[b]{\smash{{\SetFigFont{34}{40.8}{\familydefault}{\mddefault}{\updefault}{\color[rgb]{0,0,0}$b$}%
}}}}
\put(11551,-6361){\makebox(0,0)[b]{\smash{{\SetFigFont{34}{40.8}{\familydefault}{\mddefault}{\updefault}{\color[rgb]{0,0,0}$v$}%
}}}}
\put(9826,-3136){\makebox(0,0)[b]{\smash{{\SetFigFont{34}{40.8}{\familydefault}{\mddefault}{\updefault}{\color[rgb]{0,0,0}$a$}%
}}}}
\put(7426,-6436){\makebox(0,0)[b]{\smash{{\SetFigFont{34}{40.8}{\familydefault}{\mddefault}{\updefault}{\color[rgb]{0,0,0}$u$}%
}}}}
\put(9301,-5536){\makebox(0,0)[b]{\smash{{\SetFigFont{34}{40.8}{\familydefault}{\mddefault}{\updefault}{\color[rgb]{0,0,0}$z$}%
}}}}
\end{picture}%

%% file: HamVenn2a.pdf_t
\begin{picture}(0,0)%
\includegraphics{HamVenn2a.pdf}%
\end{picture}%
\setlength{\unitlength}{3947sp}%
\begingroup\makeatletter\ifx\SetFigFont\undefined%
\gdef\SetFigFont#1#2#3#4#5{%
  \reset@font\fontsize{#1}{#2pt}%
  \fontfamily{#3}\fontseries{#4}\fontshape{#5}%
  \selectfont}%
\fi\endgroup%
\begin{picture}(12645,9036)(4709,-9110)
\put(11551,-6286){\makebox(0,0)[b]{\smash{{\SetFigFont{34}{40.8}{\familydefault}{\mddefault}{\updefault}{\color[rgb]{0,0,0}$b$}%
}}}}
\put(9526,-8911){\makebox(0,0)[b]{\smash{{\SetFigFont{34}{40.8}{\familydefault}{\mddefault}{\updefault}{\color[rgb]{0,0,0}$v$}%
}}}}
\put(9826,-3136){\makebox(0,0)[b]{\smash{{\SetFigFont{34}{40.8}{\familydefault}{\mddefault}{\updefault}{\color[rgb]{0,0,0}$a$}%
}}}}
\put(7426,-6436){\makebox(0,0)[b]{\smash{{\SetFigFont{34}{40.8}{\familydefault}{\mddefault}{\updefault}{\color[rgb]{0,0,0}$u$}%
}}}}
\put(9301,-5536){\makebox(0,0)[b]{\smash{{\SetFigFont{34}{40.8}{\familydefault}{\mddefault}{\updefault}{\color[rgb]{0,0,0}$z$}%
}}}}
\end{picture}%

%% file: HamVenn2b.pdf_t
\begin{picture}(0,0)%
\includegraphics{HamVenn2b.pdf}%
\end{picture}%
\setlength{\unitlength}{3947sp}%
\begingroup\makeatletter\ifx\SetFigFont\undefined%
\gdef\SetFigFont#1#2#3#4#5{%
  \reset@font\fontsize{#1}{#2pt}%
  \fontfamily{#3}\fontseries{#4}\fontshape{#5}%
  \selectfont}%
\fi\endgroup%
\begin{picture}(12645,9036)(4709,-9110)
\put(11551,-6286){\makebox(0,0)[b]{\smash{{\SetFigFont{34}{40.8}{\familydefault}{\mddefault}{\updefault}{\color[rgb]{0,0,0}$b$}%
}}}}
\put(9526,-8911){\makebox(0,0)[b]{\smash{{\SetFigFont{34}{40.8}{\familydefault}{\mddefault}{\updefault}{\color[rgb]{0,0,0}$v$}%
}}}}
\put(9826,-3136){\makebox(0,0)[b]{\smash{{\SetFigFont{34}{40.8}{\familydefault}{\mddefault}{\updefault}{\color[rgb]{0,0,0}$a$}%
}}}}
\put(7426,-6436){\makebox(0,0)[b]{\smash{{\SetFigFont{34}{40.8}{\familydefault}{\mddefault}{\updefault}{\color[rgb]{0,0,0}$u$}%
}}}}
\put(9301,-5536){\makebox(0,0)[b]{\smash{{\SetFigFont{34}{40.8}{\familydefault}{\mddefault}{\updefault}{\color[rgb]{0,0,0}$z$}%
}}}}
\end{picture}%